\pdfoutput=1

\documentclass[a4paper,11pt,reqno]{amsart}%
\usepackage{amsfonts}
\usepackage{amsthm}
\usepackage{hyperref}
\usepackage{amssymb}
\usepackage{graphicx}
\usepackage{float}
\usepackage{cite}
\usepackage{tikz}
\usetikzlibrary{arrows}
\usepackage{xspace}
\usepackage{subcaption}
\usepackage{multirow}

\usepackage{eqparbox,array}
\usepackage{dcolumn}
\usepackage{amsmath}
\usepackage{verbatim}%
\newcolumntype{r}{D{.}{.}{-1}}

\newcommand\OSCAR{\textsc{OSCAR}\xspace}
\newcommand\Julia{\textsc{Julia}\xspace}
\newcommand\Singular{\textsc{Singular}\xspace}

\newcommand\cE{{\mathcal E}}

\newcounter{truefigure}
\makeatletter 
\renewcommand{\p@subfigure}{}
\makeatother

\newcommand{\mycaption}[1]{%
    \refstepcounter{truefigure}%
    \caption*{\textsc{Figure}~\arabic{truefigure}. #1}
}

\hypersetup{
    colorlinks,
    linkcolor={red!70!black},
    citecolor={green!60!black},
    urlcolor={blue!80!black}
}

\newtheorem{theorem}{Theorem}[section]

\theoremstyle{plain}

\theoremstyle{definition}
\newtheorem{definition}[theorem]{Definition}
\newtheorem{example}[theorem]{Example}

\theoremstyle{plain}
\newtheorem{lemma}[theorem]{Lemma}

\theoremstyle{remark}
\newtheorem{remark}[theorem]{Remark}
\numberwithin{equation}{section}

\usepackage[noend]{algorithmic}

\usepackage[linesnumbered,longend,ruled,vlined]{algorithm2e}

\usepackage[draft=true]{minted}
\setminted{breaklines=true}
\definecolor{lbcolor}{rgb}{0.9,0.9,0.9}
\setminted{bgcolor=lbcolor}
\setminted{fontsize=\footnotesize}

\usepackage{tabularx,float,capt-of}

\DeclareMathOperator{\Aut}{Aut}

\DeclareMathOperator{\mult}{mult}
\DeclareMathOperator{\trop}{trop}
\DeclareMathOperator{\coeff}{coeff}
\DeclareMathOperator{\ev}{ev}
\DeclareMathOperator{\vir}{vir}
\DeclareMathOperator{\val}{val}
\DeclareMathOperator{\lo}{loop}
\DeclareMathOperator{\Coef}{Coef}
\DeclareMathOperator{\ft}{ft}
\DeclareMathOperator{\sgn}{sgn}

\DeclareMathOperator{\preimg}{preimg}

\newcommand{\tmfloatcontents}{}
\newlength{\tmfloatwidth}
\newcommand{\tmfloat}[5]{
	\renewcommand{\tmfloatcontents}{#4}
	\setlength{\tmfloatwidth}{\widthof{\tmfloatcontents}+1in}
	\ifthenelse{\equal{#2}{small}}
	{\setlength{\tmfloatwidth}{0.45\linewidth}}
	{\setlength{\tmfloatwidth}{\linewidth}}
	\begin{minipage}[#1]{\tmfloatwidth}
		\begin{center}
			\tmfloatcontents
			\captionof{#3}{#5}
		\end{center}
\end{minipage}}

\newcommand{\cC}{\mathcal{C}}

\newcommand {\PP}{{\mathbb P}}

\newcommand {\ZZ}{{\mathbb Z}}

\newcommand {\CC}{{\mathbb C}}

\begin{document}
\title{Algorithms for Gromov-Witten Invariants of Elliptic Curves}
\author[Aga, Böhm, Hoffmann, Markwig, Traore]{Firoozeh Aga and Janko Böhm and Alain Hoffmann and Hannah Markwig and Ali Traore}
\address{Firoozeh Aga, Fraunhofer ITWM and Fachbereich Mathematik, RPTU Kaiserslautern-Landau, Postfach 3049, 67653 Kaiserslautern, Germany}
\email{aga@mathematik.uni-kl.de}
\address{Janko B\"ohm, Fachbereich Mathematik,
RPTU Kaiserslautern-Landau, Postfach 3049, 67653 Kaiserslautern, Germany}
\email{boehm@mathematik.uni-kl.de}
\address{Alain Hoffmann, Fachbereich Mathematik,
RPTU Kaiserslautern-Landau, Postfach 3049, 67653 Kaiserslautern, Germany}
\email{hoffmann@mathematik.uni-kl.de}
\address{Hannah Markwig, Fachbereich Mathematik, Eberhard Karls Universit\"at T\"ubingen}
\email{hannah@math.uni-tuebingen.de}
\address{Ali Traore, Fachbereich Mathematik, RPTU Kaiserslautern-Landau, Postfach 3049, 67653 Kaiserslautern, Germany}
\email{atraore@rptu.de}

\thanks{2010 Math subject classification: Primary 14J33, 14N35, 14T05, 81T18; Secondary 11F11, 14H30, 14N10, 14H52, 14H81}
\thanks{Gefördert durch die Deutsche Forschungsgemeinschaft (DFG) - Projektnummer 286237555 - TRR 195 [Funded by the Deutsche Forschungsgemeinschaft (DFG, German Research Foundation) - Project- ID 286237555 - TRR 195]. The work of FA was supported by Fraunhofer ITWM, Kaisers-lautern. The work of JB and HM was supported by Project A13 of SFB-TRR 195. The work of AT was supported by the German Academic Exchange Service (DAAD)  through the Mathematics in Industry and Commerce (MIC) program.}

\begin{abstract}
We present an enhanced algorithm for exploring mirror symmetry for elliptic curves through the correspondence of algebraic and tropical geometry, focusing on Gromov-Witten invariants of elliptic curves and, in particular, Hurwitz numbers. We present a new highly efficient algorithm for computing generating series for these numbers. We have  implemented the algorithm both using \Singular and \OSCAR. The implementations outperform by far the current method provided in \Singular. The \OSCAR implementation, benefiting in particular from just-in-time compilation, again by far outperforms the implementation of the new algorithm in \Singular. This advancement in computing the Gromov-Witten invariants facilitates a study of number theoretic and geometric  properties of the generating series, including quasi-modularity and homogeneity.
\end{abstract}
\maketitle

\tableofcontents

	\section{Introduction}

Mirror symmetry is a deep duality relation motivated from physics. It relates invariants of a manifold and its ''mirror manifold'' in a way allowing an exchange of methods. Tropical methods find a way into mirror symmetry most prominently in the famous Gross-Siebert programm \cite{GS06}. 

Here, we present the special case of elliptic curves, for which mirror symmetry is best understood \cite{Dij95}. The tropical side of mirror symmetry for elliptic curves was discovered within the Collaborative Research Center SFB-TRR 195  \cite{BBBM13, BGM18b, BGM18}, which is also responsible for the development of \OSCAR. 

One consequence of mirror symmetry for elliptic curves is the equality of the generating function of Hurwitz numbers (resp.\ more generally, of descendant Gromov-Witten invariants) to certain Feynman integrals which are complex analytic path integrals whose construction is governed by combinatorics. 

\emph{Hurwitz numbers} are traditional enumerative invariants that count covers of Riemann surfaces satisfying fixed ramification data. Their definition goes back to Hurwitz himself and was used to study the irreducibility of the famous moduli space of smooth genus $g$ curves, a fundamental object in algebraic geometry. In modern mathematics, Hurwitz numbers provide fruitful interactions between several mathematical areas such as geometry, topology, representation theory, combinatorics and mathematical physics \cite{CM16}. 

\emph{Tropical geometry} can be viewed as a degenerate version of algebraic geometry which has been used successfully to solve enumerative problems in algebraic geometry \cite{Mi03}. In particular, by counting \emph{tropical covers} (see Definition \ref{def-tropcover}), we can determine Hurwitz numbers via a so-called \emph{correspondence theorem} \cite{CJM10, BBBM13}. It allow to translate any computation we perform involving tropical covers to traditional Hurwitz theory.

In this chapter, we describe how one can use tropical mirror symmetry for elliptic curves in order to compute generating series for tropical Hurwitz numbers of an elliptic curve and their generalizations (up to fixed order) by means of computing Feynman integrals. We discuss explicit algorithms and implementations to handle this task. In \cite{BBBM13} an algorithm to compute Hurwitz numbers via Feynman integrals has been described which works through computing coefficients in the Laurent series expansion of the propagator of the Feynman integral in terms of the integration variables. Based on this algorithm, we develop in this chapter a significantly improved version which shows a better performance when applied to current research problems. This algorithm works by assigning to each cover a so-called flip signature by considering the permutation assigning vertices of the source tropical curve to branch points on the target elliptic curve, observing that covers with the same signature lead to the same integral. Further improvements are achieved by direct computation of coefficients of the Laurent series, avoiding explicit expansion altogether. We have created two structurally identical implementations of the algorithm, one using the computer algebra system \Singular and one using \OSCAR. Both implementations outperform the previous algorithm of \cite{BBBM13} in its implementation in \Singular by far, which demonstrates the significant structural improvement in the algorithm. The \OSCAR implementation again by far outperforms the implementation of the new algorithm in \Singular. We attribute this to the just-in-time (JIT) compilation features, which turns at run-time easy-to-comprehend \Julia code into code which executes almost as fast machine code, but also to superior arithmetic in \OSCAR. This demonstrates the potential of the \OSCAR platform. 
	
	\section{Mirror Symmetry for Elliptic Curves}
	
	\subsection{Hurwitz Numbers}
	Hurwitz numbers count branched covers of non-singular curves with a given
	ramification
	profile over fixed points. Here, we consider covers of elliptic curves. Hurwitz numbers are topological invariants,
	in particular they do not depend on the position of the branch points. 
	Moreover, since all complex elliptic curves are
	homeomorphic to the real torus, numbers of covers of an
	elliptic curve do not depend on the choice of the base curve. 
	We thus fix an arbitrary complex elliptic curve $\cE$.
	
	Let $\cC$ be a non-singular curve of genus $g$ and $\phi:\cC\rightarrow\cE$ a
	cover.
	We denote by $d$ the degree of $\phi$, i.e., the
	number of preimages of a generic point in $\cE$.
	For our purpose, it is sufficient to
	consider covers which are \emph{simply ramified}, that is, over any branch
	point exactly two sheets of the map come together and all others stay separate.
	In other words, the \emph{ramification profile} (that is, the partition of the
	degree indicating the multiplicities of the inverse images of a branch point) of
	a simple branch point is $(2,1,\ldots,1)$. 
	It follows from the Riemann-Hurwitz formula (see e.g.\ \cite{Har77}, Corollary IV.2.4) that a simply ramified cover of $\cE$ has exactly $2g-2$ branch points.
	Two covers $\varphi:\cC\rightarrow\cE$ and
	$\varphi':\cC'\rightarrow\cE$ are isomorphic if there exists an isomorphism of
	curves $\phi:\cC\rightarrow\cC'$ such that $\varphi'\circ\phi=\varphi$.

	\begin{definition}[Hurwitz numbers]\label{def-alghurwitz}
		Fix $2g-2$ points $p_1,\ldots,p_{2g-2}$ in $\cE$. We define the \emph{Hurwitz number} $N_{d,g}$ to be the weighted number of (isomorphism classes of) simply ramified covers $\phi:\cC\rightarrow \cE$ of degree $d$, where $\cC$ is a connected curve of genus $g$, and the branch points of $\phi$ are the points $p_i$, $i=1,\ldots,2g-2$. We count each such cover $\phi$ with weight $|\Aut(\phi)|$.
	\end{definition}
	%For more details on branched covers of elliptic curves and Hurwitz numbers, see e.g.\ \cite{RY10}.
	
	\begin{remark}\label{rem-branchpointsmarked}
		Note that by convention, we fix a marking of the branch points $p_i$ in this definition. In the literature, one can also find definitions which do not follow this convention and do not mark the branch points, leading to a factor of $(2g-2)!$ when compared to our definition.
	\end{remark}

	\begin{definition}
		We package the Hurwitz numbers of Definition \ref{def-alghurwitz} into a generating series as follows:
		$$ F_g(q):=\sum_{d=1}^\infty N_{d,g} q^{2d}.$$
	\end{definition}
	
	\subsection{Tropical Hurwitz Numbers}
	
	\begin{definition}[Tropical curves]
		A (generic) \emph{tropical curve} $C$ (without ends) is a connected, finite, trivalent, metric graph.
		Its
		\emph{genus} is given by its first Betti number. The
		\emph{combinatorial type} of a curve is its homeomorphism class, that is, the
		underlying graph without lengths on the edges.
	\end{definition}
	
	A tropical elliptic curve consists of one edge forming a circle of certain
	length. We will fix a tropical elliptic curve $E$ (for example, one having length~$1$).
	Moreover, to fix notation, in the following we denote by $C$ a tropical curve of genus $g$ and
	combinatorial type $\Gamma$.
	
	\begin{definition}[Tropical covers]\label{def-tropcover}
		A map $\pi:C\rightarrow E$ is a \emph{tropical cover} of $E$ if it is
		continuous,
		non-constant, integer affine on each edge and respects the \emph{balancing
			condition} at every vertex $P\in C$:
		
		For an edge $e$ of $C$ denote by $w_e$ the \emph{weight of $e$}, that is, the
		absolute
		value of the slope of $\pi_{|e}$. Consider a (small) open neighbourhood $U$ of
		$p=\pi(P)$. Then $U$ combinatorially consists of $p$ together with two rays $r_1$  and $r_2$, left and right of $p$.
		Let $V$ be the connected component of $\pi^{-1}(U)$ which contains $P$. Then
		$V$ consists of $P$ adjacent to three rays. Then $\pi$ is balanced at $P$ if
		$$\sum_{R \mbox{{ \small maps to }} r_1}
		w_R=\sum_{R \mbox{ {\small maps to }} r_2}w_R,$$ where $R$ goes over the three rays adjacent to
		$P$ and rays inherit their weights from the corresponding edges.
		
		The \emph{degree} of a cover is the weighted number of preimages of a generic
		point: For all $p\in E$ not having any vertex of $C$ as preimage we have
		$$d=\sum_{P\in C:\pi(P)=p} w_{e_P},$$ where $e_P$ is the edge of $C$ containing
		$P$. The images of the vertices of $C$ are called the \emph{branch points} of
		$\pi$.
	\end{definition}
	
\begin{figure}
		\begin{center}
			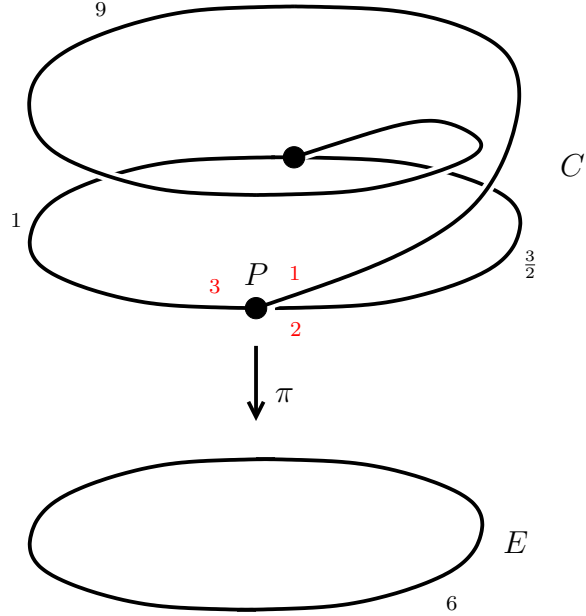
		\end{center}
		\caption[A tropical cover]{A tropical cover of degree $4$ with genus $2$ source
			curve.}\label{fig-tropcover}
	\end{figure}

	\begin{example}
		A tropical cover of degree $4$ with a genus $2$
		source curve is depicted in Figure \ref{fig-tropcover} (taken from \cite{BBBM13}). The red numbers close to the vertex $P$ are the weights of the corresponding edges, the black numbers denote the lengths. The cover is balanced at
		$P$ since there is an edge of weight $3$ leaving in one direction and an
		edge of
		weight $2$ plus an edge of weight $1$ leaving in the opposite direction.
	\end{example}
	
	We can see that the length of an edge of $C$ is determined by its weight and
	the length of its image. 
	We will therefore not specify edge lengths in the following.
	
	\begin{definition}[Isomorphisms of curves and covers]
		An \emph{isomorphism of tropical curves} is an isometry of metric graphs. Two
		covers $\pi:C\rightarrow E$ and $\pi':C'\rightarrow E$ are isomorphic if there
		is an isomorphism of curves $\phi:C\rightarrow C'$ such that $\pi'\circ\phi=\pi$.
	\end{definition}
	
	As usual when counting tropical objects we have to weight them with a certain
	multiplicity.
	
	\begin{definition}[Multiplicities]\label{def-mult}
		The multiplicity of a cover $\pi:C\rightarrow E$ is defined to be
		$$\mult(\pi):=\frac{1}{|\Aut(\pi)|}\prod_e w_e,$$ where the product goes over all edges $e$ of $C$ and $\Aut(\pi)$ is the automorphism group of $\pi$.
	\end{definition}

	\begin{definition}[Tropical Hurwitz numbers]\label{def-trophurwitz}
		Fix branch points $p_1,\ldots,p_{2g-2}$ in the tropical, elliptic curve $E$.
		The \emph{tropical Hurwitz number} $N_{d,g}^{\trop}$ is the weighted number of
		isomorphism classes of degree $d$ covers $\pi:C\rightarrow E$ having their branch
		points at the $p_i$, where $C$ is a curve of genus $g$:
		$$N_{d,g}^{\trop}:=\sum_{\pi:C\rightarrow E}\mult(\pi).$$
	\end{definition}

	\begin{theorem}[Correspondence Theorem]\label{thm-corres}
		The algebraic and tropical Hurwitz numbers of simply ramified covers of an elliptic curve coincide (see Definition \ref{def-alghurwitz} and \ref{def-trophurwitz}), i.e., we have 
		$$N_{d,g}^{\trop}=N_{d,g}.$$ 
	\end{theorem}
	
	For a proof, see the extension of the result of \cite{BBM10} in \cite{BBBM13}.

	\subsection{Feynman Integrals}
	\begin{definition}[The propagator]\label{def-prop}
		We define the \emph{propagator} $$P(z,q):=\frac{1}{4\pi^2}\wp(z,q)+\frac{1}{12}E_2(q^2)$$ in terms of the Weierstra\ss{}-P-function $\wp$ and the Eisenstein series $$E_2(q):=1-24\sum_{d=1}^\infty \sigma(d)q^d.$$ Here, $\sigma=\sigma_1$ denotes the sum-of-divisors function $\sigma(d)=\sigma_1(d)=\sum_{m|d}m$.
	\end{definition}
	The variable $q$ above should be considered as a coordinate of the moduli space of elliptic curves, the variable $z$ is the complex coordinate of a fixed elliptic curve. (More precisely, $q=e^{i\pi \tau}$, where $\tau \in \CC$ is the parameter in the upper half plane in the well-known definition of the Weierstra\ss{}-P-function.)
	
	\begin{definition}[Feynman graphs and integrals]\label{def-int}
		A \emph{Feynman graph} $\Gamma$ of genus $g$ is a trivalent connected graph of genus $g$. For a Feynman graph, we throughout fix a reference labeling $x_1,\ldots,x_{2g-2}$ of the $2g-2$ trivalent vertices and a reference labeling $q_1,\ldots,q_{3g-3}$ of the edges of $\Gamma$.
		
		For an edge $q_k$ of $\Gamma$ connecting the vertices $x_i$ and $x_j$, we define a function
		$$P_k:=P(z_i-z_j,q),$$ 
		where $P$ denotes the propagator of Definition \ref{def-prop} (the choice of sign i.e., $z_i-z_j$ or $z_j-z_i$ plays no role).
		Pick a total ordering $\Omega$ of the vertices and starting points of the form $iy_1,\ldots, iy_{2g-2}$ in the complex plane, where the $y_j$ are pairwise different small real numbers.
		We define integration paths $\gamma_1,\ldots,\gamma_{2g-2}$ by $$\gamma_j:[0,1]\rightarrow \mathbb{C}:t\mapsto iy_j+t,$$ such that the order of the real coordinates $y_j$ of the starting points of the paths equals $\Omega$. 
		We then define the integral \begin{equation}I_{\Gamma,\Omega}(q):= \int_{z_j\in \gamma_j} \prod_{k=1}^{3g-3} \left(-P_k\right).\label{eq-Igamma}\end{equation}
		Finally, we define
		$$I_\Gamma(q)= \sum_{\Omega} I_{\Gamma,\Omega}(q),$$ where the sum runs over all $(2g-2)!$ orders of the vertices. 
	\end{definition}
	
	The following relates Hurwitz numbers and Feynman integrals and can be viewed as a consequence of mirror symmetry (see Theorem 3 of \cite{Dij95}):
	\begin{theorem}[Mirror Symmetry for elliptic curves]\label{thm-mirror}
		Let $g>1$. For the definition of the invariants, see Definitions \ref{def-alghurwitz} and \ref{def-int}. We have  $$ F_g(q)=\sum_{d=1}^\infty N_{d,g} q^{2d}= \sum_{\Gamma} I_\Gamma(q)\cdot \frac{1}{|\Aut(\Gamma)|},$$ where $\Aut(\Gamma)$ denotes the automorphism group of $\Gamma$ and the sum goes over all trivalent graphs $\Gamma$ of genus $g$.
	\end{theorem}
	\begin{theorem}\label{propa}
		After coordinate change, the propagator $P(x,q)$ of Definition \ref{def-prop} with $x=e^{i \pi z}$ equals
		$$ P(x,q)=-\dfrac{x^2}{(x^2-1)^2}-\sum_{n=1}^{\infty}\left( \sum_{d\mid n} d \left(x^{2d}+x^{-2d}\right)\right) q^{2n}.$$
	\end{theorem}
 For a proof, see \cite[Theorem 2.22]{BBBM13}.

To better relate Feynman integrals with counts of tropical covers, we label the sources of our tropical covers just like a Feynman graph and consider labeled tropical covers. In that way, we can also refine the notion of degree: we can consider the degree within each edge of the source curve. Also, we can add information about the edge in the definition of the Feynman integral in Equation \ref{eq-Igamma} by making the propagator associated to edge $k$ dependent on a formal variable $q_k$, see Definition~\ref{def-Feynman} below.
 
	\begin{definition}
		We fix once and for all a base point $p_0$ in E. For a tuple $\underline{a}=(a_1,\dots,a_{3g-3})$ of non-negative integers, we define $N_{\underline{a},g}^{\trop}$ 
		to be the weighted number of labeled tropical covers  $\hat{\pi} :  C\rightarrow E$ of degree 	$\textstyle\sum_{i=1}^{3g-3} a_i$  where $C$ has genus $g$, such that $\hat{\pi}$ has its branch points at
		the prescribed positions and satisfying the condition 
		$$\# (\hat{\pi}^{-1}(p_0)\cap q_i)\cdot w_i=a_i$$
		for all $i=1,\dots,3g-3$.

	\end{definition}

	\begin{theorem}
		[Refined version of Tropical Mirror Symmetry] Let $g > 1$. We have
		$$ F_g(q_1,\dots,q_{3g-3})=\sum_{\underline{a}} N_{\underline{a},g}^{\trop}q^{2d}= \sum_{\Gamma} I_\Gamma(q_1,\dots,q_{3g-3}),$$ that is, the coefficient of the monomial $q^{2\cdot a}$ in  $I_\Gamma(q_1,\dots,q_{3g-3})$ equals 	
		$N_{\underline{a},\Gamma}^{\trop}$.
	\end{theorem}
	\begin{lemma}
		Fix a Feynman graph $\Gamma$ and an order $\Omega$ , and a tuple $(q_1, \dots , q_{3g-3} )$ as in Definition 14. We express the coefficient of $q^{2 \cdot \underline{a}}$
		in $I_{\Gamma,\Omega}(q_1,\dots,q_{3g-3})$. Assume $k$  is such that the entry $a_k = 0$, and 
		assume the edge $q_k$ connects the two vertices $x_{k_1}$
		and $x_{k_2}$. Choose the notation of the
		two vertices $x_{k_1}$ and $x_{k_2}$ such that the chosen order $\Omega$ implies  $\left| \frac{x_{k_1}}{x_{k_2}}\right|<1$  for the starting
		points on the integration paths. Then the coefficient of $q^{2 \cdot \underline{a}}$
		equals the constant term of the series
		\begin{equation}\label{eq propa}
			\prod_{k\mid a_k=0} \left(\sum_{w_k=1}^{\infty}w_k \cdot  
			\left(\frac{x_{k_1}}{x_{k_2}} \right)^{2w_k}\right)
			\cdot \prod_{k\mid a_k \neq 0} \left(\sum_{w_k \mid a_k} w_k \cdot \left( \left(\frac{x_{k_1}}{x_{k_2}} \right)^{2w_k}+ \left(\frac{x_{k_2}}{x_{k_1}} \right)^{2w_k} \right)\right)
		\end{equation}  
	\end{lemma}
	For a Taylor series  $$  F(x_1, \dots,x_n) =\sum_{\substack{1\leq j\leq n \\ a_j\geq 0}} \alpha(a_1,\dots,a_n) x_1^{a_1} ,\dots ,x_n^{a_n}
	$$
	we write $$ \coeff_{[x_1^{a_1}, \dots ,x_n^{a_n}]}(F)=\alpha(a_1,\dots,a_n).$$
	\begin{lemma} \label{coeff}
		Fix a Feynman graph $\Gamma$ and an order $\Omega$ as in Definition \ref{def-int} and write $$P_{\Gamma,\Omega} = \prod_{k=1}^{3g-3} \left(-P\left(\frac{x_{k_1}}{x_{k_2}},q\right)\right).  $$Then we have 
		
		$$ I_{\Gamma,\Omega}(q)=\coeff_{[x_1^0,\dots,x_{2g-2}^0]} (P_{\Gamma,\Omega}).$$ 
		
	\end{lemma}
	
	\subsection{Descendant Gromov--Witten Invariants}
	
	Theorem \ref{thm-mirror} can be generalized to descendant Gromov-Witten invariants, which, following Okounkov-Pandharipande \cite{OP06}, are related to covers with more complicated ramification profiles.
	
	A  {\it stable map} of degree $d$ from a curve of genus $g$  to $\cE$ with $n$ markings is a map $f: \cC \to \cE$, where  $\cC$ is a connected projective curve with at worst nodal singularities, and with $n$ distinct nonsingular marked points $x_1,\ldots,x_n\in \cC$, such that $f_\ast([\cC]) = d[\cE]$ and $f$ has a finite group of automorphism. 
	The moduli space of stable maps, denoted  $\overline{\mathcal{M}}_{g,n} (\cE,d)$, is a proper Deligne-Mumford stack of virtual dimension $2g-2+n$ \cite{Beh97,BF97}. 
	The $i$-th  evaluation morphism is the map $\ev_i: \overline{\mathcal{M}}_{g,n} (\cE,d) \to \cE$ that sends a point $[\cC, x_1, \ldots, x_n, f]$ to  $f(x_i)
	\in \cE$. 
	The $i$-th cotangent line bundle $\mathbb{L}_i \to \overline{\mathcal{M}}_{g,n} (\cE,d)$ is defined by a canonical identification of its fiber over a moduli point $(\cC, x_1, \ldots, x_n, f)$
	with the cotangent space $T^\ast_{x_i}(\cC)$.  The first Chern class of the cotangent line bundle is called a {\it psi class} ($\psi_i= c_1(\mathbb{L}_i)$).
	
	\begin{definition}
		Fix $g,n,d$ and let $k_1, \ldots, k_n$ be non-negative integers with $$k_1+\ldots+k_n = 2g-2.$$ The {\it stationary descendant Gromov-Witten invariant}  $\langle \tau_{k_1}(pt) \ldots \tau_{k_n}(pt) \rangle_{g,n}^{\cE,d}$ is  defined by
		\begin{equation}
			\langle \tau_{k_1}(pt) \ldots \tau_{k_n}(pt) \rangle_{g,n}^{\cE,d} = \int_{[\overline{\mathcal{M}}_{g,n} (\cE,d)]^{\vir}} \prod_{i=1}^n ev_i^\ast(pt) \psi_i^{k_i},
		\end{equation}
		where $pt$ denotes the class of a point in $\cE$.
	\end{definition}
	
	In order to define tropical multiplicities, we also need to discuss relative descendant Gromov-Witten invariants of $\mathbb{P}^1$.
	They are constructed using moduli spaces of {\it relative stable maps} $\overline{\mathcal{M}}_{g,n} (\PP^1, \mu,\nu,d)$, where part of the data specified are the ramification profiles $\mu$ and $\nu$ which we fix over $0$ resp.\ $\infty\in \PP^1$. The preimages of $0$ and $\infty$ are marked. A detailed discussion of spaces of relative stable maps to $\PP^1$ and their boundary is not necessary for our purpose, we refer to \cite{Vak08}. We use operator notation and denote
	\begin{equation}
		\langle \mu| \tau_{k_1}(pt) \ldots \tau_{k_n}(pt)|\nu \rangle_{g, n}^{\PP^1,d} = \int_{[\overline{\mathcal{M}}_{g,n} (\PP^1, \mu,\nu,d)]^{\vir}} \prod_{i=1}^n \ev_i^\ast(pt) \psi_i^{k_i}.
	\end{equation}

	Tropically, adding descendants amounts to allowing vertices of higher genus and higher valency. For this purpose, we first generalize the definition of abstract tropical curve:

 \begin{definition}
	An \emph{abstract tropical} \emph{curve} is a connected metric graph $\Gamma$, such that edges leading to leaves (called \emph{ends}) have infinite length, together with a genus function $g:\Gamma\rightarrow \ZZ_{\geq 0}$ with finite support. Locally around a point $p$, $\Gamma$ is homeomorphic to a star with $r$ halfrays. 
	The number $r$ is called the \emph{valence} of the point $p$ and denoted by $\val(p)$. We identify the vertex set of $\Gamma$ as the points where the genus function is nonzero, together with points of valence different from $2$. The vertices of valence greater than $1$ are called  \textit{inner vertices}. Besides \emph{edges}, we introduce the notion of \emph{flags} of $\Gamma$. A flag is a pair $(V,e)$ of a vertex $V$ and an edge $e$ incident to it ($V\in \partial e$). 
	Edges that are not ends are required to have finite length and are referred to as \emph{bounded} or \textit{internal} edges.
	
	A \emph{marked tropical curve} is a tropical curve whose leaves are labeled. An isomorphism of a tropical curve is an isometry respecting the leaf markings and the genus function. The \emph{genus} of a tropical curve $\Gamma$ is given by
	\[
	g(\Gamma) = h_1(\Gamma)+\sum_{p\in \Gamma} g(p)
	\]
	The \emph{combinatorial type} is the equivalence class of tropical curves obtained by identifying any two tropical curves which differ only by edge lengths.
 \end{definition}
	
	Definition \ref{def-tropcover} introducing tropical covers easily extends to allow such more general source curves.
		
	\begin{definition}[Psi- and point conditions]\label{def-PsiAndPointConditions}
		We say that a tropical cover $\pi:\Gamma_1\rightarrow \Gamma_2$ with a marked end $i$ {\it satisfies a psi-condition} with power $k$ at $i$, if the vertex $V$ to which the marked end $i$ is adjacent has valency  $k+3-2g(V)$. We say $\pi:\Gamma_1\rightarrow \Gamma_2$ {\it satisfies the point conditions} $p_1,\dots,p_n\in\Gamma_2$ if $$\lbrace\pi(1),\dots,\pi(n)\rbrace=\lbrace p_1,\dots,p_n \rbrace.$$
	\end{definition}
	
	Fix $g,n,d$ and let $k_1, \ldots, k_n$ be non-negative integers with $$k_1+\ldots+k_n = 2g-2.$$
	Let $\pi:\Gamma\rightarrow E$ be a tropical cover such that $\Gamma$ is of genus $g$ and has $n$ marked ends. 
	Fix $n$ distinct points $p_1,\ldots,p_n\in E$.
	Assume that at the marked end~$i$, a psi-condition with power $k_i$ is satisfied, and that the point conditions are satisfied. 
	The marked ends must be adjacent to different vertices, since they satisfy different point conditions. It follows from an Euler characteristic argument incorporating the valencies imposed by the psi-conditions that $\Gamma$ has exactly $n$ vertices, each adjacent to one marked end.

	Locally at the marked end $i$, the cover sends the vertex to an interval consisting of two flags $f$ and $f'$. We define the {\it local vertex multiplicity} $\mult_i(\pi)$ to be a one-point relative descendant Gromov-Witten invariant: 
	\begin{equation} \mult_i(\pi)= \langle \mu_f| \tau_{k_i}(pt)|\mu_{f'} \rangle_{g_i, 1}^{\PP^1,d_i},\label{eq-localmult}\end{equation}
	where $g_i$ denotes the genus of the vertex adjacent to the marked end $i$, $d_i$ its local degree, and $\mu_f$ resp. $\mu_{f'}$ the ramification profiles above the two flags of the image interval.
	
	We define the multiplicity of $\pi$ to be
	\begin{equation}
		\frac{1}{|\Aut(\pi)|}\cdot \prod_i \mult_i(\pi)\cdot \prod_e \omega(e).\label{eq-mult}
	\end{equation}
	
	Note that all ends of a tropical cover of $E$ are contracted ends, with image points the points $p_i$ we fix as conditions in $E$.
	
	\begin{definition}[Tropical stationary descendant Gromov-Witten invariant of $E$]
		For $g$, $n$, $d$, $k_1,\ldots,k_n$ as above, define the \emph{tropical stationary descendant Gromov-Witten invariant }
		$$\langle \tau_{k_1}(pt) \ldots \tau_{k_n}(pt) \rangle_{g,n}^{\cE,d,\trop}$$
		to be the weighted count of tropical genus $g$ degree $d$ covers of $E$ with $n$ marked points satisfying point and psi-conditions as above, each counted with its multiplicity as defined in (\ref{eq-mult}). 
		
	\end{definition}
	
	Also for descendant invariants, a correspondence theorem holds:
	\begin{theorem}[Correspondence Theorem for descendants]\label{thm-corres-desc}
		A stationary descendant Gromov-Witten invariant of $\cE$ coincides with its tropical counterpart:
		$$\langle \tau_{k_1}(pt) \ldots \tau_{k_n}(pt) \rangle_{g,n}^{\cE,d}=\langle \tau_{k_1}(pt) \ldots \tau_{k_n}(pt) \rangle_{g,n}^{E,d,\trop}.$$
	\end{theorem}
	For a proof, see \cite[Theorem 3.2.1]{CJMR16}. 
	
	Also generating series of descendant invariants can be expressed in terms of Feynman integrals. However, we as well need to extend our notion of Feynman integrals to accommodate this setting.
	
	\begin{definition}[The propagator and the $\mathcal{S}$-function]\label{def-props}
		We define the \emph{propagator} as a (formal) series in $x$ and $q$,
		$$
		P(x,q)= \sum_{w=1}^\infty w\cdot x^{w} +\sum_{a=1}^\infty \left(\sum_{w|a}w \left(x^{w}+ x^{-w}\right)\right)q^{a}
		$$
		and the $\mathcal{S}$-function as series in $z$,
		$$\mathcal{S}(z)=\frac{\sinh(z/2)}{z/2}.$$
		
		We also introduce a further formal power series in $q$, which should be viewed as the propagator for loop edges:
		
		$$P^{\lo}(q)= \sum_{a=1}^\infty \Biggl(\sum_{w|a} w\Biggl) q^a.
		$$
	\end{definition}
	
	Analogous to the way we generalized our notion of abstract tropical curves which can serve as sources of our tropical covers, we have to generalize the notion of Feynman graphs. We allow any graph $\Gamma$ without ends with $n$ vertices which are labeled $x_1,\ldots,x_n$ and with labeled edges $q_1,\ldots,q_r$. By convention, we assume that $q_1,\ldots,q_s$ are loop edges and $q_{s+1},\ldots,q_r$ are non-loop edges.

	\begin{definition}[Feynman integrals] \label{def-Feynman} Let $\Gamma$ be a Feynman graph.
		Let $\Omega$ be an order of the $n$ vertices of $\Gamma$.
		For $k>s$, denote the vertices adjacent to the (non-loop) edge $q_k$ by $x_{k^1}$ and $x_{k^2}$, where we assume $x_{k^1}<x_{k^2}$ in $\Omega$.
		We define the \emph{Feynman integral} for $\Gamma$ and $\Omega$ to be
		$$I_{\Gamma,\Omega}(q)=\Coef_{[x_1^{0}\ldots x_n^{0}]} \prod_{k=1}^s P^{\lo}(q) \cdot  \prod_{k=s+1}^{r} P\Big(\frac{x_{k^1}}{x_{k^2}},q\Big).$$
		%and the \emph{refined Feynman integral} to be
		$$I_{\Gamma,\Omega}(q_1,\ldots,q_r)=\Coef_{[x_1^{0}\ldots x_n^{0}]}\prod_{k=1}^s P^{\lo}(q_k)   \prod_{k=s+1}^{r} P\Big(\frac{x_{k^1}}{x_{k^2}},q_k\Big).$$
		Finally, we set 
		$$I_{\Gamma}(q)=\sum_\Omega I_{\Gamma,\Omega}(q),$$
		where the sum goes over all $n!$ orders of the vertices of $\Gamma$.
		
	\end{definition}
	If we assume $|x|<1$ to express the (in $q$) constant coefficient of the (non-loop) propagator (that is, the first sum appearing in the propagator series in Definition~\ref{def-prop}) as the rational function $\frac{x^2}{\left( x^2-1\right)^2}$ (using geometric series expansion), we can view the  series from which we take the $x_1^{0}\ldots x_n^{0}$-coefficient in the Feynman integral above as a function on a Cartesian product of elliptic curves. The Feynman integral then becomes a path integral in complex analysis as above.
	Note that using the change of coordinates $x=e^{i\pi u}$ the (non-loop) propagator has the form of Definition \ref{def-prop}.
	
	\begin{definition}[Feynman integrals with vertex contributions]\label{def-Feynmanvertex}
		Let $\Gamma$ be a Feynman graph, and equip it with an additional genus function $\underline{g}$ associating a nonnegative integer $g_i$ to every vertex $x_i$. Let $\Omega$ be an order of the $n$ vertices of $\Gamma$.
		We adapt our notion of propagators from Definitions \ref{def-props} and \ref{def-Feynman} to include vertex contributions:
		for non-loop edges, we set
		\begin{align*} \tilde{P}(\frac{x_{k^1}}{x_{k^2}},q) =&
			\sum_{w=1}^\infty \mathcal{S}(w z_{k^1} )\mathcal{S}(w z_{k^2}) \cdot w \cdot  \left(\frac{x_{k^1}}{x_{k^2}}\right)^w \\& + \sum_{a=1}^\infty \Bigg(\sum_{w|a}
			\mathcal{S}(w z_{k^1} )\mathcal{S}(w z_{k^2})\cdot w \cdot \left(\left(\frac{x_{k^1}}{x_{k^2}}\right)^w+ \left(\frac{x_{k^2}}{x_{k^1}}\right)^w\right)\Bigg)\cdot q^a.\end{align*}
		For loop-edges connecting the vertex $x_{k^1}$ to itself, we set
		$$\tilde{P}^{\lo}(q)= \sum_{a=1}^\infty\Bigg(\sum_{w|a}\mathcal{S}(w z_{k^1} )^2\cdot w\Bigg) q^a.
		$$
		The variables $z_{k_i}$ are new variables introduced for each vertex in order to take care of the genus contribution.
		
		We define the \emph{Feynman integral with vertex contributions} for $\Gamma$, $\underline{g}$ and $\Omega$ to be 
\begin{align*}I_{\Gamma,\underline{g},\Omega}(q)= \Coef_{[z_1^{2g_1}\ldots z_n^{2g_n}]} \Coef_{[x_1^{0}\ldots x_n^{0}]} \prod_{i=1}^n\frac{ 1}{\mathcal{S}(z_i)} \prod_{k=1}^s  \tilde{P}^{\lo}(q) \prod_{k=s+1}^{r} \tilde{ P}(\frac{x_{k^1}}{x_{k^2}},q).\end{align*}

		Again, we set
		$$I_{\Gamma,\underline{g}}(q)=\sum_\Omega I_{\Gamma,\underline{g},\Omega}(q) $$
		where the sum goes over all $n!$ orders of the vertices.
		%and
		%$$I_{\Gamma,\underline{g}}(q_1,\ldots,q_r)= \sum_\Omega I_{\Gamma,\underline{g},\Omega}(q_1,\ldots,q_r).$$
		
	\end{definition}

	With this, we can generalize Theorem \ref{thm-mirror} to the version involving descendants:

	\begin{theorem}[Mirror symmetry for $E$, version with descendants]\label{thm-mirror-desc}
		Fix $g\geq 2$, $n\geq 1$ and $k_1, \ldots, k_n\geq 1$ satisfying $k_1+\ldots+k_n = 2g-2$.
		
		We can express the series of descendant Gromov-Witten invariants of $E$ in terms of Feynman integrals,
		
		$$\sum_{d\geq 1}\langle \tau_{k_1}(pt) \ldots \tau_{k_n}(pt) \rangle_{g,n}^{\cE,d}q^d = \sum_{(\ft(\Gamma),\underline{g})} \frac{1}{|\Aut(\ft(\Gamma),{\underline{g}})|}I_{\Gamma,\underline{g}}(q),$$
		where $\Gamma$ is a Feynman graph with a genus function $\underline{g}$, such that the vertex $x_i$ has genus $g_i$ and valency $k_i+2-2g_i$, and such that $h^1(\Gamma)+\sum g_i=g$, and where we consider automorphisms of unlabeled graphs ($\ft$ is the forgetful map that forgets all labels of a Feynman graph $\Gamma$) that are required to respect the genus function.
	\end{theorem}
	
	For a proof, see \cite{Li11, BGM18}.
	
	\section{Computation of Generating Series}\label{algbasic}
    
	In this section, we will discuss our improved algorithm for computing tropical Hurwitz numbers and Gromov-Witten invariants via evaluation of Feynman integrals. The algorithm implemented in the library \texttt{ellipticcovers.lib} \cite{BBBM} for the computer algebra system \Singular, is based on evaluating the propagator product as a Laurent series in the variables $x_i$ in the correct order given by $\Omega$. For $k>s$, denoting the vertices adjacent to the (non-loop) edges $q_k$ by $x_{k^1}$ and $x_{k^2}$ in any order, and fixed degree partition $\underline{a}$ over the base point, we evaluate 		$$I_{\Gamma,\Omega,\underline{a}}=\Coef_{[x_{\Omega(1)}^{0}\ldots x_{\Omega(n)}^{0}]}\Coef_{\underline{a}}\prod_{k=1}^s P^{\lo}(q_k)   \prod_{k=s+1}^{r} P\Big(\frac{x_{k^1}}{x_{k^2}},q_k\Big)$$
  extracting the constant coefficient in the $x_i$ and sum over all $\Omega$. Vertex contributions are introduced as in Definition \ref{def-Feynmanvertex}. For more details, see also \cite{BBBM13}.

	\subsection{Basic Improved Algorithm}\label{alghurwitz}
	Before turning to the flip signature, we first describe a structural improvement to the above mentioned method.
	Let us consider an edge labeled $q_k$ in the graph $\Gamma$, connecting vertices $x_i$ and $x_j$. The propagator associated with this edge will be a polynomial in the variable $q_k$. However, we are specifically interested in extracting the coefficient of the term $q_k^{2a_k}$, where $a_k$ is the $k$-th entry of $\underline{a}$. Write   $N=\sum_{i=1}^{3g-3} a_i$.
	By Lemma  \ref{coeff} and Equation (\ref{eq propa}),  the coefficient of $q^{2 \cdot \underline{a}}$ in $I_{\Gamma,\Omega}(q_1,\dots,q_{3g-3})$
	can explicitely described as follows. For $a_k > 0$, we have
		
		$$\sum_{w_k \mid a_k} w_k \cdot \left( \left( \frac{x_{k_1}}{x_{k_2}}
		\right)^{w_k} + \left( \frac{x_{k_2}}{x_{k_1}} \right)^{w_k} \right) =
		\frac{\underset{w_k \mid a_k}{\sum}  w_k \cdot (x_{k_1}^{N + w_k} x_{k_2}^{N
				- w_k} + x_{k_1}^{N - w_k} x_{k_2}^{N + w_k})}{x_{k_1}^N x_{k_2}^N}$$		
	and for $a_k = 0$, we obtain		
		$$\underset{w_k = 1}{\overset{N}{\sum}} w_k \cdot \left(
		\frac{x_{k_1}}{x_{k_2}} \right)^{w_k} = \frac{\underset{w_k =
				1}{\overset{N}{\sum}} w_k \cdot x_{k_1}^{N + w_k} x_{k_2}^{N -
				w_k}}{x_{k_1}^N x_{k_2}^N}.$$
	In our computations, the focus lies on computing the numerator of this representation of the propagator. The  denominator shift  $ x_i^N x_j^N$ can be taken into account separately within the function \texttt{coefficient\_of\_term} which determines the coefficient of a specific monomial $q^{2 \cdot \underline{a}}$ by extracting the coefficient of the $ x_i^N x_j^N$ term of the numerator. 
	\subsection{Generating Series for Hurwitz Numbers}
	\subsubsection{Flip Signature}
 In our approach, we introduce a flip signature with respect to the vertex ordering $\Omega$ to group covers having the same Feynman integral $I_{\Gamma,\Omega}$ together into a single computation. This grouping provides a further significant speedup of the computation of generating series of Gromov-Witten invariants Feynman integral $I_\Gamma = \sum_\Omega I_{\Gamma,\Omega}$ by identifying how many permutations $\Omega$ lead to the same signature, and determining the integral for one representative for each signature.
 
             The permutation $\Omega$ determines which vertices of the edges will occur in the numerator and which of them will occur in the denominator in the corresponding propagator. While there are $(2g-2)!$ orders of branch points, there are at most $2^{2g-2}$ possible flip signatures, two possibilities for each edge of the graph. Moreover, due to symmetries in the  propagator product typically the number of flip signatures is even smaller.
            
            We define the flip signature as follows: We first find the $\emph{sign}$ of an edge $q_k$ for a given vertex ordering $\Omega$.
            We set the sign of an edge $q_k$ to be 1, if $\Omega(i) < \Omega(j)$, else it is set to be -1.
            For a given graph $\Gamma$, order $\Omega$, and branch type  $\underline{a} = (a_1, \ldots ,a_{3g-3} )$,  Algorithm \ref{alg:sgn} \texttt{flip\_signature} determines the flip signature as a vector of length $3g-3$ assigning each edge the value $-1$ in case the degree of the edge is zero and its sign is $-1$ or the value $0$ otherwise. Note that for edges of positive degree permutation of $x_i$ and $x_k$ does not influence the integral since the above specified propagator summand is symmetric in $x_i$ and $x_k$, hence it does not depend on the choice of order of Laurent series expansion in $x_i$ and $x_k$.
            
            For a given graph $\Gamma$, order $\Omega$, and branch type  $\underline{a} = (a_1, \ldots ,a_{3g-3} )$,  Algorithm \ref{alg:sgn} \texttt{flip\_signature} determines the flip signature as a vector of length $3g-3$ assigning each non-loop edges the value $-1$ in case the degree of the edge is zero and its sign is $-1$, and the value $0$ otherwise. To detect loop edges we assign in this case the value $-2$.
            
        \begin{algorithm}
		\caption{\texttt{flip\_signature}}\label{alg:sgn}
		\KwIn{Graph $\Gamma$, permutation $\Omega$, branch type $\underline{a}$.}
		\KwOut{Flip signature $b$.}
		\Begin{
			E = edges of $\Gamma$  ($x_{k_1}$ -- $e_k$ -- $x_{k_2}$)
			
			\ForEach{$e_{k} \text{ in } E$}
			{
			\ForEach{$a_{k}  \text{ in } a$}
				{
					\If{ $a_{k}=0$}
						{
								\If{$\preimg(\Omega,x_{k_1})<  \preimg(\Omega,x_{k_2})$}
									{
										$b_k = -1 $\
									}
								\Else{$b_k = 0$ }
						}
					\If{$x_{k_1}=x_{k_2}$}
						{
							$b_k = -2$
						} 
			    }
		}
	return $b$
	}
	\end{algorithm}
\subsubsection{Multiplicities}

In the following Algorithm \ref{alg:flip} \texttt{signature\_and\_multiplicities}, we find how often each flip signature occurs over all $\Omega$. This is then called the \emph{multiplicity} of the flip signature. Note that we only permute the vertices of edges, where a flip can occur. For the other vertices, the order will not change the flip signature. Hence, we count over a subgroup of the permutation group of all vertices, and therefore need to adjust the multiplicities of the flip signatures by the corresponding index.

\begin{algorithm}
	\caption{\texttt{signature\_and\_multiplicities}}\label{alg:flip}
	\KwIn{Graph $\Gamma$, branch type $\underline{a}$.}
	\KwOut{Vector with flip signatures and their multiplicity}
	\Begin{
		$E =$ edges of $\Gamma$  ($x_{k_1}--e_k--x_{k_2}$)
		
		$V =$ set of all vertices of edges $k$ of $\Gamma$ with $a_k = 0$
		
		$P =S(V)$ all permutations of $V$

        $D = $ empty dictionary
		
		\ForEach{$\Omega \in P$}
		{
			$y = \sgn(\Gamma,\Omega,a)$
         
            \If{$y \text{ key of } D$} 
              {
               $D[y] = D[y] +1$
              } 
            \Else
             {
               $D[y] = 1$
             } 
		}
		
		multiply each value of $D$ by $\frac{(\#\text{ vertices of } \Gamma)!}{|P|}$
		
		return $D$
	}
\end{algorithm}

\subsubsection{Feynman Integral}

Building on signatures and multiplicities, we now proceed to the computation of the Feynman integral for a fixed branch type $\underline{a}$, which is one coefficient in the multivariate generating series and is denoted by  $I_{\Gamma,\underline{a}}$. We also give an algorithm to compute the multivariate generating series up to a given degree $d$, which  is denoted by  $I_{\Gamma,d}$. The respective algorithms are Algorithm \ref{alg:specific} (\texttt{feynman\_integral\_branchtype}) and Algorithm \ref{alg:feynman} (\texttt{feynman\_integral\_degree}).

The functions $\operatorname{const\_term}$ and $\operatorname{non\_const\_term}$ compute the $q$-constant term of the propagator multiplied by $x_i^n x_j^n$ and the coefficient of $q_k^{2d}$ of the propagator multiplied by $x_i^n x_j^n$, respectively. For given $n$ they are given as 
$$\operatorname{const\_term}(x_i , x_j , N) = \sum_{w=1}^n w \cdot x_i^{n+w}x_j^{n-w}$$
and 
$$\operatorname{non\_const\_term}(x_i , x_j , q_k , d, N) = \sum_{w|d} w \cdot (x_i^{n+w} x_j^{n-w} + x_i^{n-w} x_j^{n+w} ) \cdot q_k^{2d}.$$ Here we can choose $n = \sum_{k=1}^{3g -3} a_k$, since any variable $x_i^d$ with degree  $d < -n$ or $n < d$ can not contribute to the constant term $x_1^0 \cdot \ldots \cdot x_{2g-2}^0$.

\begin{algorithm}
	\caption{\texttt{feynman\_integral\_branchtype}}\label{alg:specific}
	\KwIn{Graph $\Gamma$, branch type $\underline{a}$. %; $l$ leak vector (optional)
    }
	\KwOut{$I_{\Gamma,\underline{a}}$}
	\Begin{
		$E =$ edges of $\Gamma$
		
		$N = \sum_{i=1}^{3g-3} a_i $
		
		$D = \operatorname{signatures\_and\_multiplicities}(\Gamma,\underline{a})$ 
		
		%$v(D)=Value$(D)=multiplicities$
		
		%Key($D$)=flip signatures.

        $p=0$
		
		\ForEach{$f \in \operatorname{keys}(D)$}{
            $tmp = 1$

            \For{$j$ \textbf{\textup{from}} $1$ \textbf{\textup{to}} $\operatorname{size}(f)$ }
			{
                
				\If{$k_j=-1$}{ $ tmp = tmp\cdot \operatorname{const\_term}(x_{E_{j,1}}, x_{E_{j,2}}, N)$}
				\ElseIf {$k_j=0$}
					{ $ tmp = tmp\cdot \operatorname{const\_term}(x_{E_{j,2}}, x_{E_{j,1}}, N)$}
				
				\Else{$ tmp = tmp\cdot \operatorname{non\_const\_term}(x_{E_{j,1}}, x_{E_{j,2}},q_j,f_j, N)$}
				
			}	
			$p = p+D(f)\cdot  \Coef_{[x_1^{N} \cdots x_{2g-2}^{n}]} (tmp) $
		}
		return $p$
	}
\end{algorithm}

\begin{algorithm}
	\caption{\texttt{feynman\_integral\_degree}}\label{alg:feynman}
	\KwIn{Graph $\Gamma$, maximum degree $d$.}
	\KwOut{$I_{\Gamma,d}$}
	\Begin{

        $sum = 0$
       \For{$j$ \textbf{\textup{from}} $1$ \textbf{\textup{to}} $d$ }
          {
			
		$A =$ partitions of $d$ into $\#$(edges of $\Gamma$) summands.
 
		\ForEach{$a \in A$}
		{
           sum=sum$+\operatorname{feynman\_integral\_branchtype}(\Gamma,a)$
		}
        }
		return sum
	}
\end{algorithm}

	\subsection{Generating Series for Descendant Gromov-Witten Invariants}

The improved algorithms from Section \ref{alghurwitz} have a straight-forward generalization to the case of Gromov-Witten invariants with Psi-classes, following the formula given at the beginning of Section  \ref{algbasic}. We omit the respective algorithms here, but give example and timings at \cite{GWpackage}.

\section{Timings}
	In this section, we evaluate the performance of our implementation for various graphs and degrees $d$ comparing: the implementation of the original algorithm in the \Singular library \texttt{ellipticcovers.lib} in the function \texttt{gromovWitten} (denoted by Singular-1 in the tables), the implementation of the new algorithm in the function \texttt{feynman\_integral} included in an updated version of \texttt{ellipticcovers.lib} (denoted by Singular-2 in the tables), and the respective \OSCAR-based implementation in the function \texttt{feynman\_integral} of the package \texttt{GromovWitten.jl} \cite{GWpackage} (denoted by OSCAR in the tables). This package is available at \begin{center}\url{https://github.com/singular-gpispace/GromovWitten}.\end{center} All computations have been done sequentially on an M-core with $3.2$ MHz. We have verified that all answers are consistent. The tables show all data collected in the respective configuration and are in seconds rounded to $3$ digits.

\begin{figure}[ht]
    \centering
    \begin{subfigure}[b]{.48\linewidth}
        \centering
        \raisebox{-0.5\height}{\includegraphics[width=4cm]{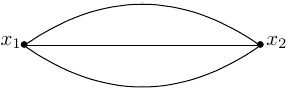}}
        \mycaption{Caterpillar genus $2$.}
        \label{fig:cat2}
    \end{subfigure}
    \hfill
    \begin{subfigure}[b]{.48\linewidth}
        \centering
        \raisebox{-0.5\height}{\includegraphics[width=3.4cm]{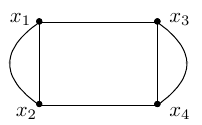}}
        \mycaption{Caterpillar genus $3$.}
        \label{fig:cat3}
    \end{subfigure}
\end{figure}

\begin{figure}[ht]
    \centering
    \begin{subfigure}[b]{.48\linewidth}
        \centering
        \raisebox{-0.5\height}{\includegraphics[width=3.4cm]{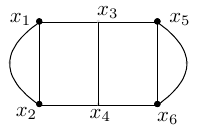}}
        \mycaption{Caterpillar genus $4$.}
        \label{fig:cat4}
    \end{subfigure}
    \hfill
    \begin{subfigure}[b]{.48\linewidth}
        \centering
        \raisebox{-0.5\height}{\includegraphics[width=4cm]{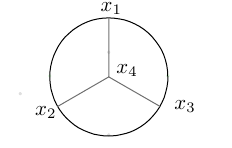}}
        \mycaption{Star graph $K_{1, 3}$.}
        \label{fig:star}
    \end{subfigure}
\end{figure}

 We consider the source curves depicted in  Figures~\ref{fig:cat2}--\ref{fig:star}. To generate the graphs considered in the timings, the following code can be used:

\begin{minted}{jlcon}
julia> using GromovWitten

# caterpillar genus 2
julia> G =  graph([(1,2), (1,2),(1,2)]);
	
# caterpillar genus 3
julia> G = graph([(1,2),(1,2),(2,4),(1,3),(3,4),(3,4)]);
	
# caterpillar genus 4
julia> G = graph((1,2),(1,2),(1,3),(2,4),(3,4),(3,5),(4,6),(5,6),(5,6));
	
# star graph K_{1,3}
julia> G = graph((1,2),(1,3),(1,4),(2,4),(2,3),(3,4));
\end{minted}  %

\begin{figure}[ht]
	\raisebox{0.0\height}{\includegraphics[width=6cm,height=5cm]{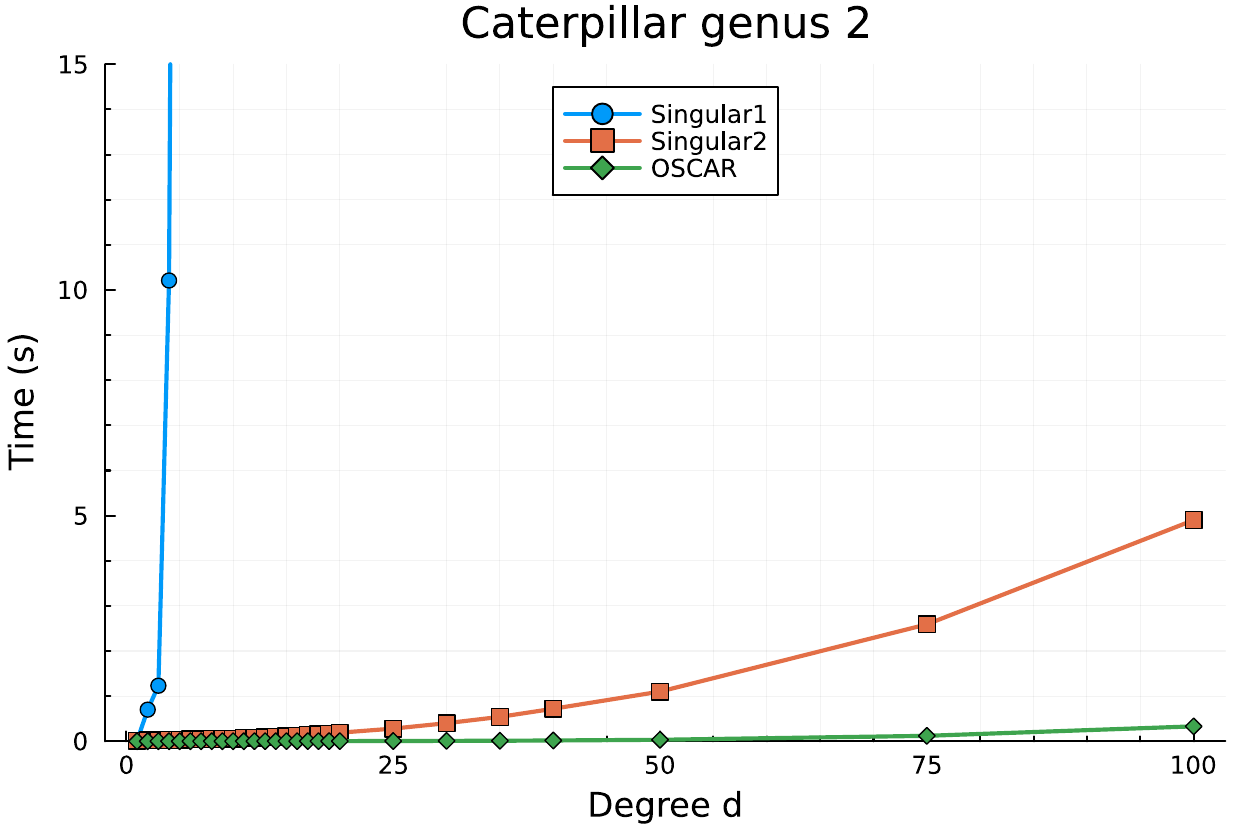}}\qquad\raisebox{0.0\height}{\includegraphics[width=6cm,height=5cm]{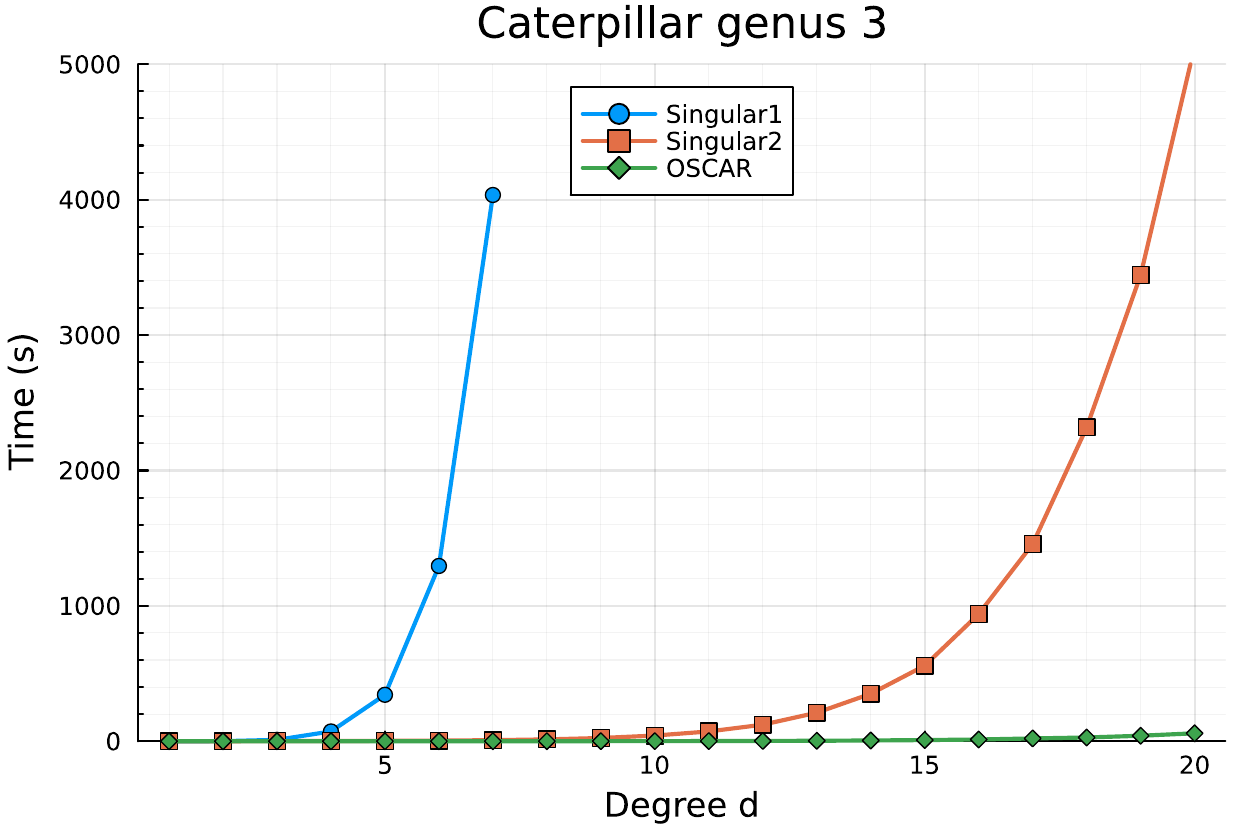}}
	
	\
	
	\raisebox{0.0\height}{\includegraphics[width=6cm,height=5cm]{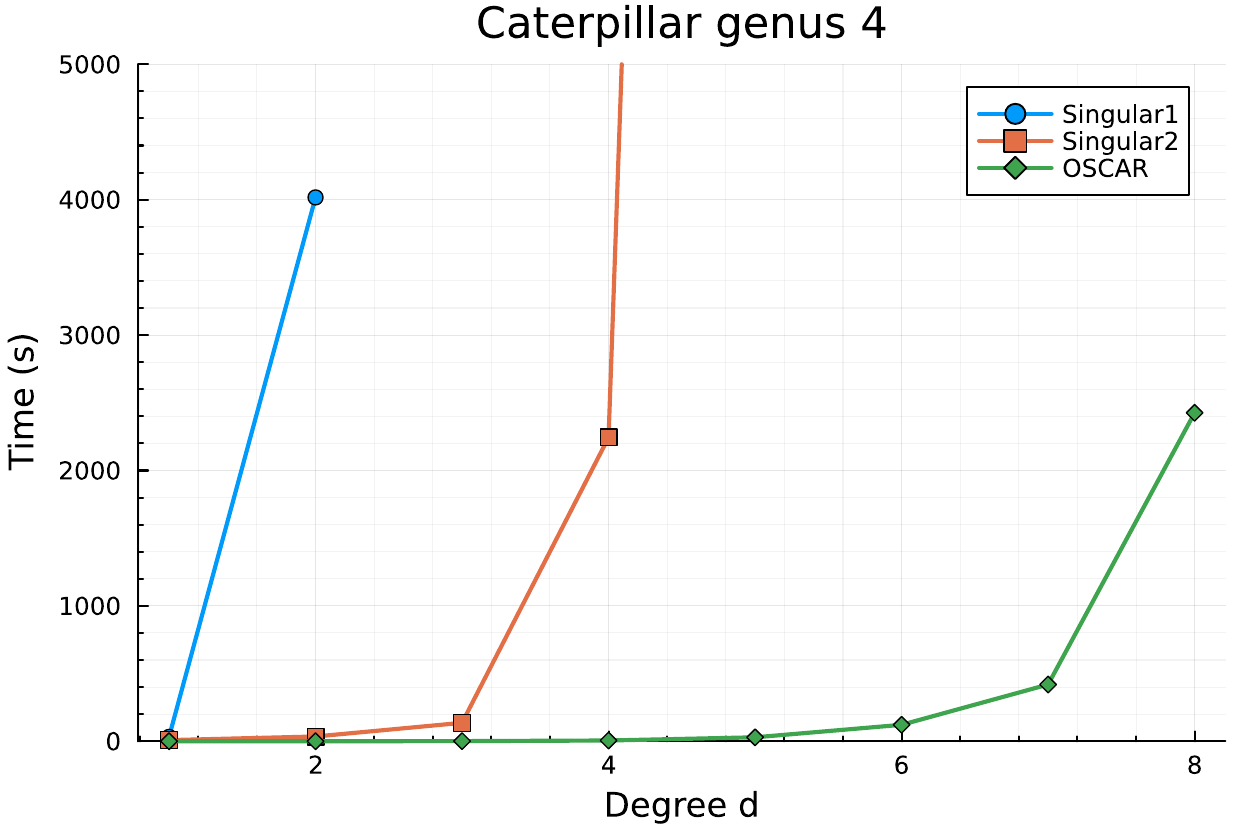}}\qquad\raisebox{0.0\height}{\includegraphics[width=6cm,height=5cm]{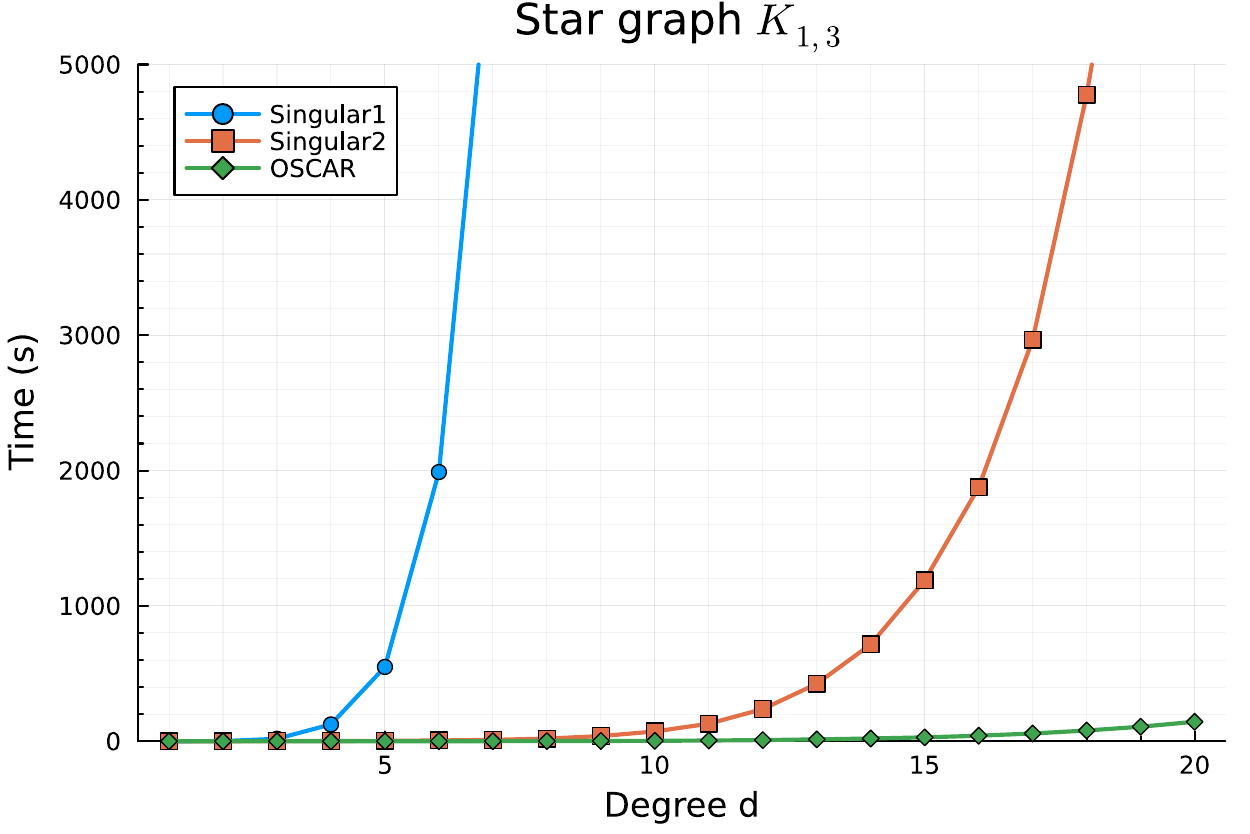}}
	\caption{Visualization of the timings for Singular-1, Singular-2, and~OSCAR.}
    \label{fig timings GW}
\end{figure}

\begin{figure}[ht]
	\raisebox{0.0\height}{\includegraphics[width=6cm,height=5cm]{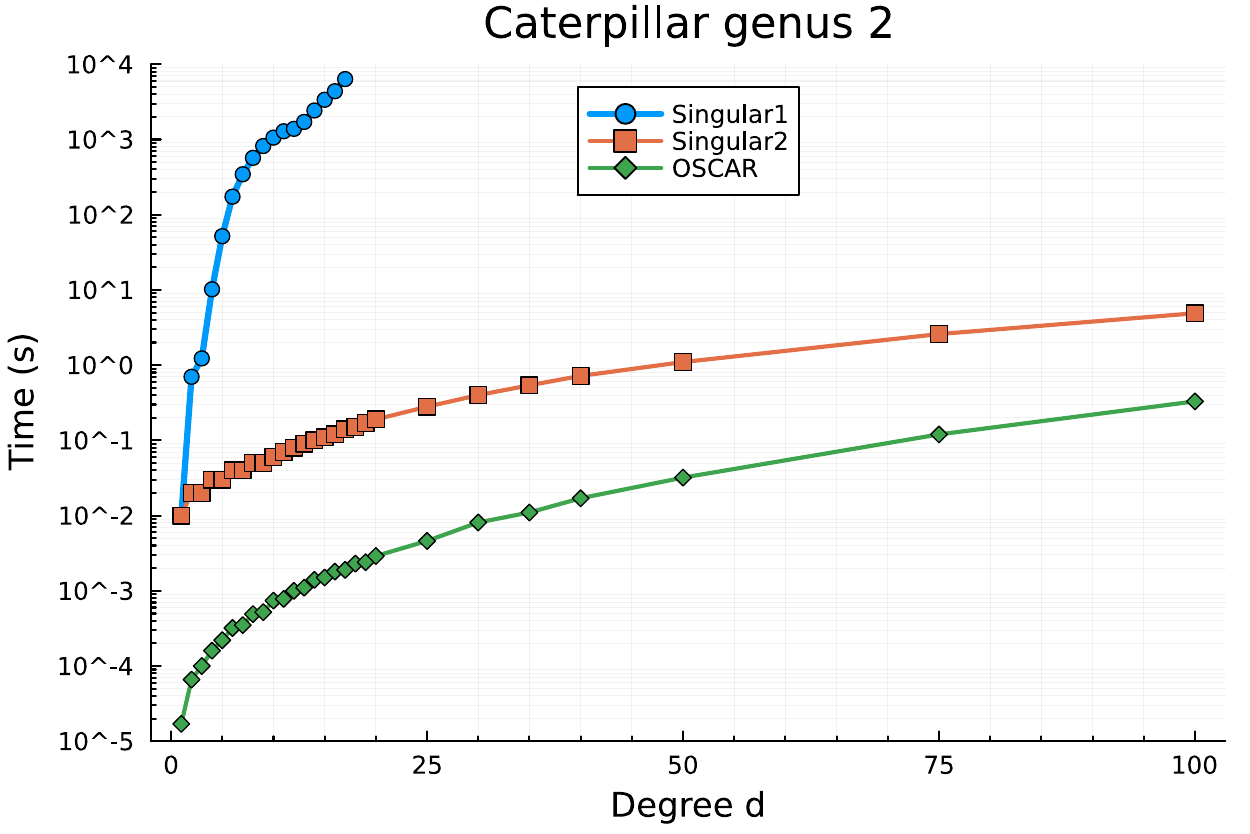}}\qquad\raisebox{0.0\height}{\includegraphics[width=6cm,height=5cm]{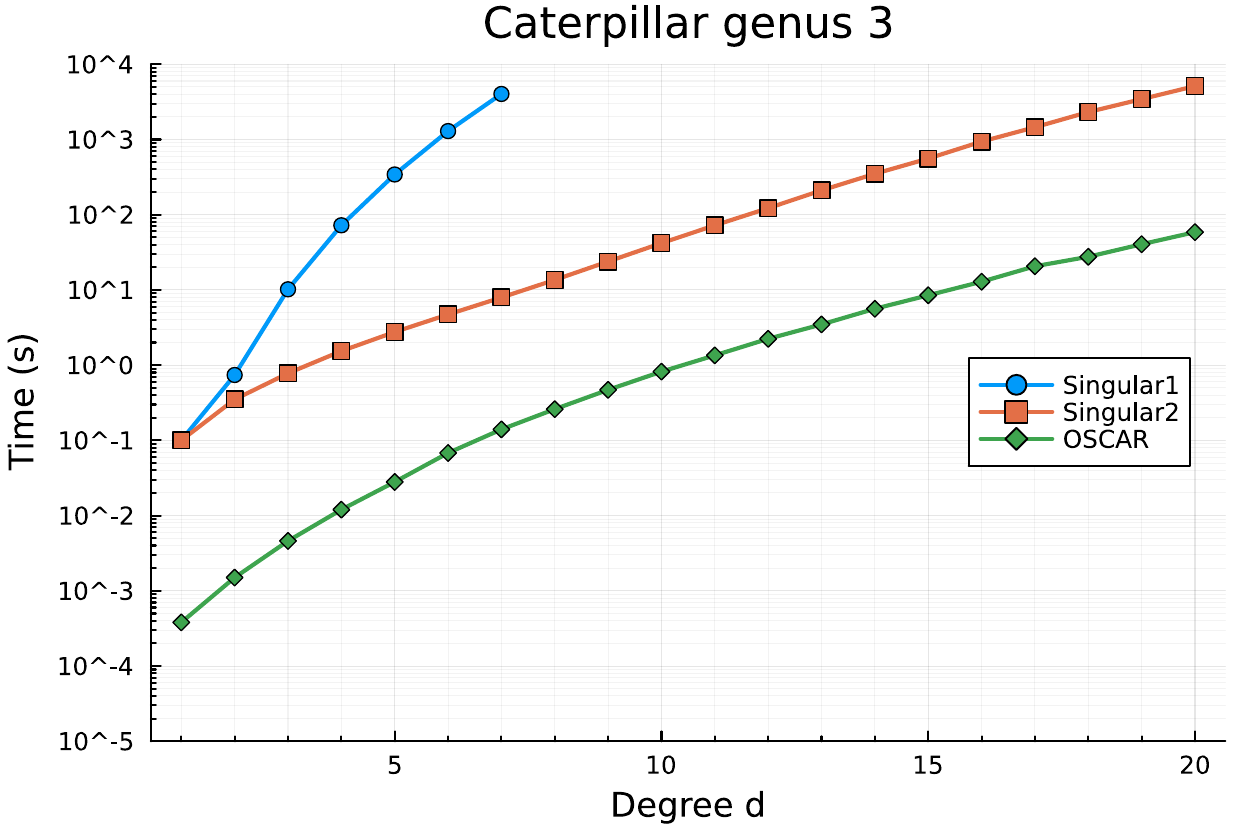}}
	
	\
	
	\raisebox{0.0\height}{\includegraphics[width=6cm,height=5cm]{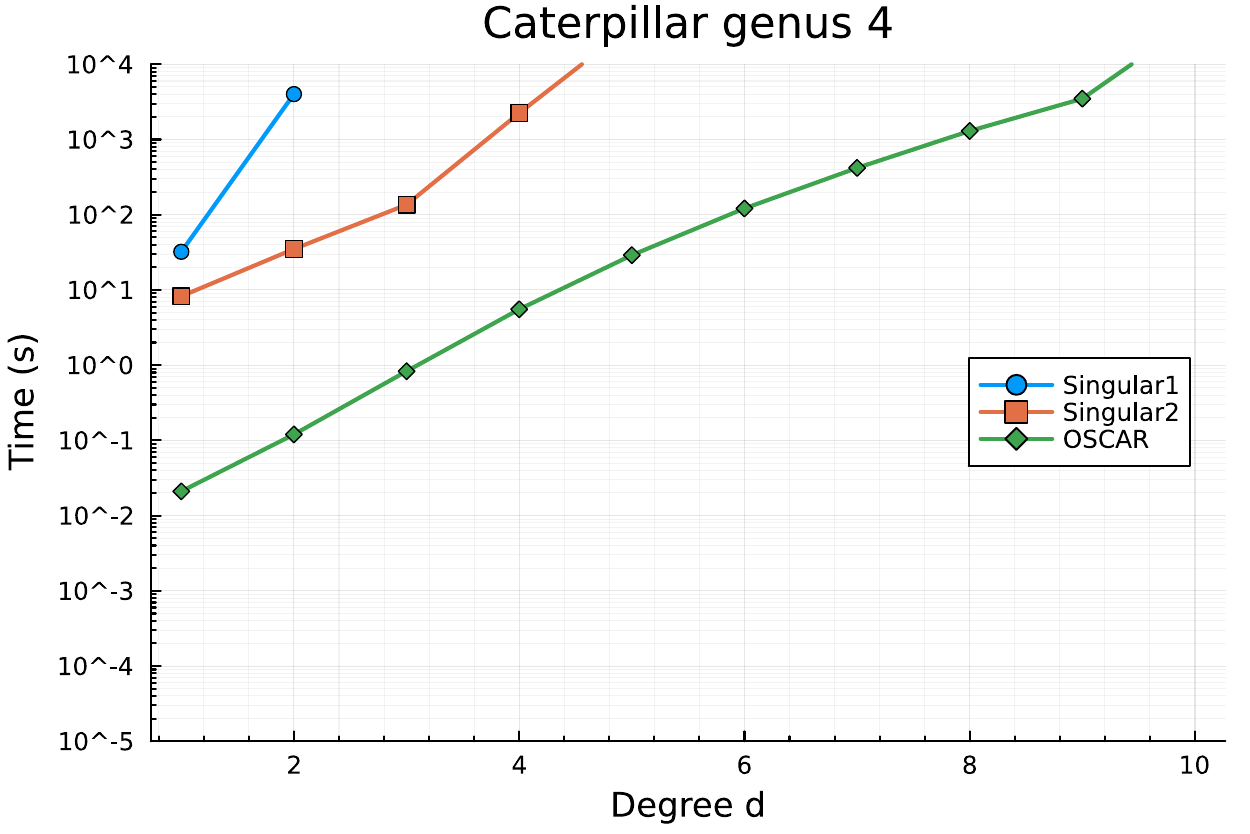}}\qquad\raisebox{0.0\height}{\includegraphics[width=6cm,height=5cm]{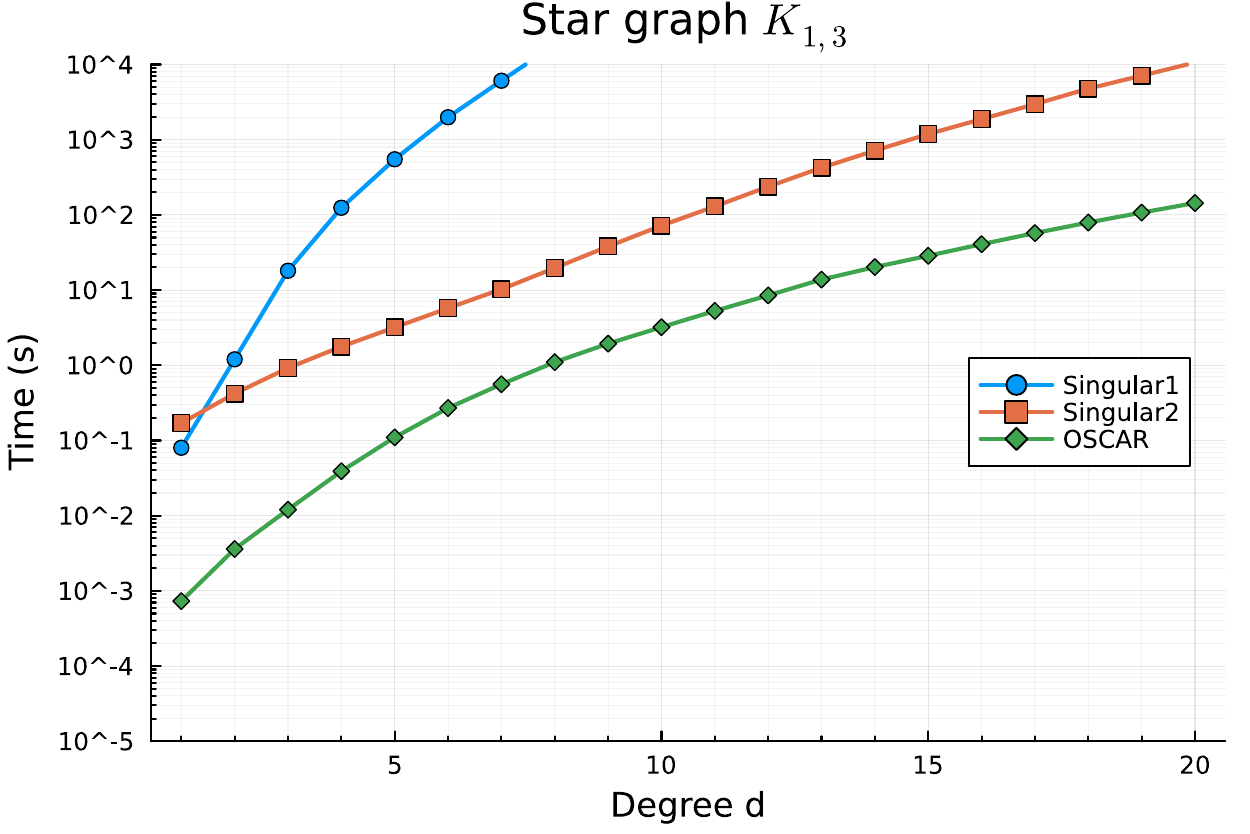}}
	\caption{Visualization of the timings for Singular-1, Singular-2, and~OSCAR (logarithmic time scale).}
    \label{fig timings GW2}
\end{figure}

For the timings, see Tables \ref{table gw1} and \ref{table gw3}, and for a visualization of the timings, see Figures~\ref{fig timings GW} and ~\ref{fig timings GW2}. We observe that the new algorithm presented in this paper outperforms the previous one by far, and that the \OSCAR implementation of that algorithm outperforms the \Singular implementation again by far, although it is structurally identical. The demonstrates the potential of the new \OSCAR platform, in particular, the benefit of just-in-time compilation. The \OSCAR implementation allows to compute until about double or more of the degree in the expansion of the formal generating series than the \Singular implementation, which again can achieve about three times the degree compared to the old algorithm. The six-fold or more increase in degree allows for example to study homogeneity properties of quasimodular representations of the generating series  for many graphs of larger genus which were previously not tractable. There is ongoing work on parallelization of the computation.

\begin{table}[ht]
	\centering
	\begin{tabularx}{\textwidth}{|c|*{6}{X|}}
		\cline{1-7}
		\multirow{2}{*}{degree d} &  \multicolumn{3}{c|}{Caterpillar genus 2} & \multicolumn{3}{c|}{Caterpillar genus 3} \\
		\cline{2-7}
		& \mbox{Singular-1} & \mbox{Singular-2} & \mbox{OSCAR} & \mbox{Singular-1} & \mbox{Singular-2} & \mbox{OSCAR} \\
		\hline
  1 & 0.01 & 0.01 & 0.00001& 0.10 & 0.10 & 0.0003\\
    2 & 0.70 & 0.02 & 0.00006 & 0.74 & 0.35  & 0.001\\
    3 & 1.23& 0.021 & 0.0001 & 10.20 & 0.78  & 0.004\\
    4 & 10.21 & 0.030 & 0.00016 & 72.41 & 1.54 & 0.012\\
    5 & 51.84 & 0.032 & 0.0002 & 343.03  & 2.76  & 0.03\\
    6 & 173.94 & 0.041 & 0.0003 & 1300  & 4.72  & 0.68 \\
    7 & 345.16& 0.044 & 0.00035 & 4030  & 7.94  & 0.14\\
    8 & 570.06 & 0.051 & 0.0004 &--& 13.6 & 0.26\\
    9 & 819.31 & 0.053 & 0.0005 &--& 23.9 & 0.47 \\
    10 &1060& 0.06 & 0.0007 &--& 41.9  & 0.82\\
    11 &1287& 0.07 & 0.0007 &--& 72.3  &  1.35\\
    12 &1386& 0.08 & 0.0010 &--& 122  & 2.24 \\
    13 &1714& 0.09 & 0.0011 &--& 211  & 3.48 \\
    14 &2433& 0.10 & 0.0014 &--& 351 & 5.62 \\
    15 &3376& 0.11 & 0.0015 &--& 559 & 8.49\\
    16 &4390& 0.12 & 0.0018 &--& 940 & 12.9\\
    17 &--& 0.14 & 0.0019 &--& 1 450 & 20.7\\
    18 &--& 0.15 & 0.0023 &--& 2 300 &27.6 \\
    19 &--& 0.17 & 0.0024 &--& 3 440 & 40.5 \\
    20 &--& 0.19 & 0.0029 &--& 5 130 & 58.7\\

		\hline
	\end{tabularx}
	\caption{Timings for caterpillar source curves of genus 2 and genus~3.}
    \label{table gw1}
\end{table}

\begin{table}[ht]
	\centering
	\begin{tabularx}{\textwidth}{|c|*{6}{X|}}
		\cline{1-7}
		\multirow{2}{*}{degree d} &  \multicolumn{3}{c|}{Caterpillar genus 4} & \multicolumn{3}{c|}{Star graph $K_{1,3}$} \\
		\cline{2-7}
		& \mbox{Singular-1} & \mbox{Singular-2} & \mbox{OSCAR} & \mbox{Singular-1} & \mbox{Singular-2} & \mbox{OSCAR} \\
		\hline
	  1 & 32.16 & 8.24  & 0.021 &  0.08 & 0.17 & 0.0007   \\
	  2 & 4016 & 35.0 & 0.12& 1.20 & 0.42 & 0.003  \\
	  3 &--& 136 & 0.83& 18.1 & 0.92 &  0.012   \\
	  4 &--& 2240 & 5.53& 124 & 1.77 &  0.039 \\
	  5 &--& 32 890 & 29.0& 549& 3.19 & 0.112  \\
	  6 &--&--& 121& 1 990 & 5.80 & 0.27 \\
	  7 &--&--& 419& 6 080 & 10.3 & 0.56  \\
	  8 &--&--&  2 420& 18 200 & 19.7 & 1.10\\
	  9 &--&--&3 490&--& 38.0 &1.94 \\
	  10 &--&--&38 700&--& 72.0 & 3.21 \\
	  11 &--&--&--&--& 130 & 5.28  \\
	  12 &--&--&--&--& 238 &  8.50 \\
	  13 &--&--&--&--& 424 &  13.8 \\
	  14 &--&--&--&--& 716 &  20.2  \\
	  15 &--&--&--&--& 1 190  & 28.6 \\
	  16 &--&--&--&--& 1 880& 40.8  \\
	  17 &--&--&--&--& 2 970  & 57.3  \\
	  18 &--&--&--&--& 4 777  & 79.2 \\
	  19 &--&--&--&--& 7 103  & 107 \\
	  20 &--&--&--&--& 10 620  & 144 \\
	\hline
\end{tabularx}
 \caption{Timings for caterpillar genus 4 and $K_{1,3}$ source curves.}
\label{table gw3}
\end{table}

\pagebreak[4]

 \subsection{Example Computation}
 In this section, we provide example code for computing the Hurwitz nunmers for a caterpillar genus $3$ source curve.
\begin{minted}{jlcon}
julia> using GromovWitten

julia> G = FeynmanGraph([(1, 3), (1, 2), (1, 2), (2, 4), (3, 4), (3,4)] )
graph([(1, 3), (1, 2), (1, 2), (2, 4), (3, 4), (3, 4)])

julia> F = FeynmanIntegral(G)

julia> a = [0, 2, 1, 0, 0, 1];

julia> feynman_integral_branch_type(F, a)
256*q[2]^4*q[3]^2*q[6]^2

julia> feynman_integral_degree(F, 3)
288*q[1]^3 + 32*q[1]^2*q[2] + 32*q[1]^2*q[3] + 32*q[1]^2*q[5] + 32*q[1]^2*q[6] + 8*q[1]*q[2]*q[5] + 8*q[1]*q[2]*q[6] + 8*q[1]*q[3]*q[5] + 8*q[1]*q[3]*q[6] + 24*q[2]^3 + 152*q[2]^2*q[3] + 8*q[2]^2*q[5] + 8*q[2]^2*q[6] + 152*q[2]*q[3]^2 + 32*q[2]*q[3]*q[5] + 32*q[2]*q[3]*q[6] + 32*q[2]*q[4]^2 + 8*q[2]*q[4]*q[5] + 8*q[2]*q[4]*q[6] + 8*q[2]*q[5]^2 + 32*q[2]*q[5]*q[6] + 8*q[2]*q[6]^2 + 24*q[3]^3 + 8*q[3]^2*q[5] + 8*q[3]^2*q[6] + 32*q[3]*q[4]^2 + 8*q[3]*q[4]*q[5] + 8*q[3]*q[4]*q[6] + 8*q[3]*q[5]^2 + 32*q[3]*q[5]*q[6] + 8*q[3]*q[6]^2 + 288*q[4]^3 + 32*q[4]^2*q[5] + 32*q[4]^2*q[6] + 24*q[5]^3 + 152*q[5]^2*q[6] + 152*q[5]*q[6]^2 + 24*q[6]^3

julia> f = feynman_integral_degree_sum(F, 3)
288*q[1]^6 + 32*q[1]^4*q[2]^2 + 32*q[1]^4*q[3]^2 + 32*q[1]^4*q[5]^2 + 32*q[1]^4*q[6]^2 + 8*q[1]^4 + 8*q[1]^2*q[2]^2*q[5]^2 + 8*q[1]^2*q[2]^2*q[6]^2 + 8*q[1]^2*q[3]^2*q[5]^2 + 8*q[1]^2*q[3]^2*q[6]^2 + 24*q[2]^6 + 152*q[2]^4*q[3]^2 + 8*q[2]^4*q[5]^2 + 8*q[2]^4*q[6]^2 + 152*q[2]^2*q[3]^4 + 32*q[2]^2*q[3]^2*q[5]^2 + 32*q[2]^2*q[3]^2*q[6]^2 + 8*q[2]^2*q[3]^2 + 32*q[2]^2*q[4]^4 + 8*q[2]^2*q[4]^2*q[5]^2 + 8*q[2]^2*q[4]^2*q[6]^2 + 8*q[2]^2*q[5]^4 + 32*q[2]^2*q[5]^2*q[6]^2 + 8*q[2]^2*q[6]^4 + 24*q[3]^6 + 8*q[3]^4*q[5]^2 + 8*q[3]^4*q[6]^2 + 32*q[3]^2*q[4]^4 + 8*q[3]^2*q[4]^2*q[5]^2 + 8*q[3]^2*q[4]^2*q[6]^2 + 8*q[3]^2*q[5]^4 + 32*q[3]^2*q[5]^2*q[6]^2 + 8*q[3]^2*q[6]^4 + 288*q[4]^6 + 32*q[4]^4*q[5]^2 + 32*q[4]^4*q[6]^2 + 8*q[4]^4 + 24*q[5]^6 + 152*q[5]^4*q[6]^2 + 152*q[5]^2*q[6]^4 + 8*q[5]^2*q[6]^2 + 24*q[6]^6

julia> substitute(f)
1792*q[1]^6 + 32*q[1]^4 
\end{minted}  %

\bibliographystyle {plain}
\bibliography {bibliographie}

\end{document}